\documentclass[12pt]{amsart}
\usepackage{amsmath, amssymb}

\theoremstyle{plain}
\usepackage{amsmath}
\usepackage{amssymb}
\usepackage{amsfonts}
\newtheorem{theorem}{Theorem}

\newtheorem{lemma}{Lemma}
\newtheorem{remark}{Remark}
\newtheorem{definition}{Definition}
\newtheorem{corollary}{Corollary}
\numberwithin{equation}{section}

\def \Z {\mathbb{Z}}

\def \P {\mathbb{P}}

\def \DD {\mathcal{D}}

\def \a {\alpha}
\def \b {\beta}

\def \e {\varepsilon}

\def \d {\delta}

\def \z {\zeta}

\def \< {\langle}
\def \> {\rangle}
\def \^ {\widehat}
\def \sign {{\rm sign}}
\def \dist {{\rm dist}}

\def \Comp {{\mathit{Comp}}}
\def \Incomp {{\mathit{Incomp}}}

\newcommand{\bbA}{{\bf A}}

\newcommand{\bbH}{{\bf H}}

\newcommand{\bbI}{{\bf I}}

\newcommand{\bbW}{{\bf W}}

\newcommand{\bbX}{{\bf X}}

\newcommand{\bbx}{{\bf x}}

\newcommand{\bbY}{{\bf Y}}

\newcommand{\bbb}{{\bf b}}

\newcommand{\beq}{\begin{equation}}
\newcommand{\eeq}{\end{equation}}
\newcommand{\bqa}{\begin{eqnarray}}
\newcommand{\eqa}{\end{eqnarray}}
\newcommand{\bqn}{\begin{eqnarray*}}
\newcommand{\eqn}{\end{eqnarray*}}
\newcommand{\non}{\nonumber \\}
\newcommand{\bdes}{\begin{description}}
\newcommand{\edes}{\end{description}}

\def\underwiggle 1{\ifmmode\setbox\TempBox=\hbox{$ 1$}\else\setbox\TempBox=\hbox{1}\fi
\setbox\TempBoxA=\hbox to \wd\TempBox{\hss\char'176\hss}
\rlap{\copy\TempBox}\smash{\lower9pt\hbox{\copy\TempBoxA}} }
\begin{document}

\title[Circular law]{Circular law, Extreme Singular values and Potential theory}

\author{Guangming Pan
  \and Wang Zhou}

\thanks{
  W. Zhou. was supported by a grant
 R-155-050-055-133/101 at the National University of Singapore}

\address{Eurandom, P.O.Box 513, 5600MB Eindhoven, the Netherlands}
\email{stapgm@gmail.com}

\address{Department of Statistics and Applied Probability, National University of
 Singapore, Singapore 117546}
\email{stazw@nus.edu.sg}

\maketitle
\begin{abstract}

Consider the empirical spectral distribution of complex random $n\times n$
matrix whose entries are independent and identically distributed random
variables with mean zero and variance $1/n$. In this paper, via applying potential
theory in the complex plane and analyzing extreme singular values, we prove that
 this distribution converges, with probability
one, to the uniform distribution over the unit disk in the complex
plane, i.e. the well known circular law, under the finite fourth moment
assumption on  matrix elements.

\end{abstract}

\section{Introduction}

Let \{$X_{kj}$\}, $k,j=\cdots,$ be a double array of independent and
identically distributed (i.i.d.) complex random variables (r.v.'s) with
$EX_{11}=0$ and $ E|X_{11}|^2=1$.  The complex eigenvalues of the
matrix $n^{-1/2}\bbX=n^{-1/2}(X_{kj})$ are denoted by $\lambda_1,\cdots,\lambda_n$.
The two-dimensional empirical spectral distribution $\mu_n(x,y)$ is
defined as
\begin{equation}\label{intr1}
\mu_n(x,y)=\frac{1}{n}\sum\limits_{k=1}^nI(Re(\lambda_k)\leq x, \ Im
(\lambda_k)\leq y).
\end{equation}
The study of $\mu_n(x,y)$ is related to understanding the random behavior
of slow neutron resonances in nuclear physics. See \cite{me}. Since 1950's
it has been conjectured that, under the unit varaince condition,
$\mu_n(x,y)$ converges to the so-called {\it circular law}, i.e.  the
uniform distribution over the unit disk in the complex plane. Up to now,
this conjecture is only proved in some partial cases.

The first answer for complex normal matrices was given in \cite{me}
based on the joint density function of the eigenvalues of
$n^{-1/2}\bbX$.
 Huang in \cite{hu} reported
that this result was obtained in an unpublished paper of Silverstein
(1984). After more than one decade, Edelman \cite{edel} also showed
that the expected empirical spectral
distribution converges to the circular law for real normal
matrices. It is Girko who investigated the circular law for general
matrix with independent entries for the first time in \cite{girk1}.
But Girko imposed, not only moment conditions, but also strong smooth
conditions on matrix entries.
Later on, he further published
a series of papers (for example, \cite{girk2}) about this problem.
However, as pointed out in \cite{bai2} and \cite{gt}, Girko's argument
includes serious mathematical gaps.
The rigorous argument of the
conjecture was given by Bai in his 1997 celebrated paper
 \cite{bai2} for general random matrices.
In addition to the finite $(4+\varepsilon)th$ moment condition Bai still assumed
that  the joint density of the real and imaginary part of the entries is
bounded. Again, the result was further improved by Bai and Silverstein
under the assumption $E|X_{11}|^{2+\eta}<\infty$ in their comprehensive
book \cite{bai3}, but the finiteness condition of the density of matrix
entries is still there. Recently, G\"otze and Tikhomirov \cite{gt} gave a
proof of the convergence of $E\mu_n(x,y)$ to the circular law under the
strong moment assumption that the entries have sub-Gaussian tails or are
sparsely non-zero instead of the condition about the density of the
entries in \cite{bai2}.

Generally speaking,
there are five approaches to studying the spectral distribution of random matrices.
The difficulty of the circular conjecture is that the methodologies used in
Hermitian matrices do not work well in non-Hermitian ones. There was no
powerful tool to attack this conjecture.

1. {\it Moment method.} Moments are very important characteristics of r.v.'s.
They have many applications in probability and statistics. For example,
we have moment estimators in statistics. As far as we know,
it is Wigner \cite{wig1} \cite{wig2} who introduced moment method into random matrices.
Since then, the moment method has been very successful in
establishing the convergence of the empirical spectral distribution of
Hermitian matrices. Bai did a lot of important work. One can refer to \cite{bai3}.
But moment method fails to work in non-Hermitian ones,
because for any complex r.v. $Z$ uniformly distributed
over any disk centered at $0$, one can verify that for any $m\geq 1$
$$
EZ^m=0.
$$

2. {\it Stietjes transform.} Another powerful tool in random matrices theory is the
Stieltjes transform, which is defined by
\begin{equation}\label{intr2}
m_G(z):=\int\frac{1}{\lambda-z}dG(\lambda),\quad z\in
{\Bbb{C}}^+\equiv\{z\in {\Bbb C},\ Im (z)>0\},
\end{equation}
for any distribution function $G(x)$. The basic property of Stieltjes transform is that
it is a representing class of probability measures. This property offers one
a strong analytic machine.
 Still see \cite{bai3} and
the references therein.  However, the Stieltjes
transform of $n^{-1/2}\bbX$ is unbounded if $z$ coincides
with one eigenvalue. So this leads to serious difficulties when dealing with
the Stieltjes transform of $n^{-1/2}\bbX$.

3. {\it Orthogonal polynomials.} The study of orthogonal polynomials
goes back as far as Hermite. For the deep connections between orthogonal polynomials
and random matrices, one can refer to \cite{dei}. Orthogonal polynomials are
usually limited to Guassian random matrices.
Moreover, orthogonal polynomials are only suitable to deriving the spacing
between consecutive eigenvalues for large classes of random matrices
(see \cite{deift}).

4. {\it Characteristic functions.} There is a long history of characteristic
functions. In 1810, Laplace used Fourier transform, i.e. characteristic
functions to prove central limit theorem for bounded r.v.'s. Then
in 1934 P. L\'evy reproved Linderberg central limit theorem by characteristic
functions. From that time on, characteristic functions are well known to
almost every mathematician. Surprisingly,  one can not see any application
of characteristic functions in random matrices until 1984.   Girko combined together
the characteristic function of $\mu_n(x,y)$ and the
Stieltjes transform, trying  to prove the conjecture
in \cite{girk1}. Developing ideas proposed by Girko \cite{girk1},
 Bai reduced the conjecture to estimating the smallest singular value of
$n^{-1/2}\bbX-z\bbI$ in \cite{bai2}. However, one should note
that some uniform estimate of the smallest singular values of
$n^{-1/2}\bbX-z\bbI$ with respect to $z$ will be required if
the method in \cite{bai2} is employed.

5. {\it Potential theory.} Potential theory is the terminology given to the wide
area of analysis encompassing such topics as harmonic and subharmonic functions,
the boundary problem, harmonic measure, Green's function, potentials and capacity.
Since Doob's famous book \cite{doob} appeared, it is widely accepted that potential
theory and probability theory are closely related. For example, superharmonic
functions correspond to supermartingales.

The logarithmic potential of
a measure $\mu$ (see \cite{saff}) is defined by
\begin{equation}\label{intr3}
U^{\mu}(z):=\int \log\frac{1}{|z-t|}d\mu(t),
\end{equation}
where $\mu(t)$ is any positive finite Borel measure with support in
a compact subset of the complex plane. There is also an inversion
formula, i.e. $\mu$ can be defined through $U^{\mu}$ as $d\mu=-(2\pi)^{-1}\Delta U^{\mu}$,
where $\Delta$ is the two dimensional Laplacian operator. This relation makes
Khoruzhenko in \cite{kho} suggest to use potential theory to derive the circular law.
Then G\"otze and Tikhomirov in \cite{gt} used the
logarithmic potential of $E\mu_n$ convoluted by a smooth distribution
 to provide a proof for the convergence of $E\mu_n$ to the circular
law with entries being sub-Gaussian or sparsely non-zero.

In this paper, the conjecture, the convergence of $\mu_n(x,y)$ to the
circular law with probability one, is established under the assumption
that the underlying r.v.'s have finite fourth moment.
Compared with \cite{gt}, we work on the logarithmic potential of
$\mu_n(x,y)$ directly, while \cite{gt} depends on the logarithmic
potential of a convolution of $E\mu_n(x,y)$ and the uniform
distribution on the disk of radius $r$. 

The main result of this paper is formulated as follows.

\begin{theorem}\label{theo3}
Suppose that $\{X_{jk}\}$ are i.i.d. complex r.v.'s with
$EX_{11}=0,$ $ E|X_{11}|^2=1$ and $E|X_{11}|^4<\infty$. Then, with
probability one, the empirical spectral distribution function
$\mu_n(x,y)$ converges to the uniform distribution over the unit
disk in two dimensional space.
\end{theorem}

\begin{remark}
The bounded density condition in \cite{bai2} and
the sub-Gaussian assumption in \cite{gt} are not needed any more.
\end{remark}

Theorem \ref{theo3} will be handled by potential theory in conjunction with
estimates for the
smallest singular value of $n^{-1/2}\bbX-z\bbI$.

The research of the smallest singular values originates from von Neumann
and his colleagues.
They guessed that
\begin{equation}\label{d1}
s_n(\bbX)\sim n^{-1/2}\quad \text{with high probability},
\end{equation}
with $s_n(\bbX)$ being the smallest singular value of $\bbX$. Edelman in
\cite{ede} proved it for random Gaussian matrices, i.e., for each $\varepsilon\geq 0$
\begin{equation}\label{d2}
P(s_n(\bbX)\leq \varepsilon n^{-1/2})\sim\varepsilon.
\end{equation}
Rudelson and Vershynin in \cite{rud2}
solved it for real random matrices, i.e., for every $\delta>0$ there exist
$\varepsilon>0$ and $n_0$ depending only on
$\delta$ and the fourth moment of $X_{jk}$ so that
\begin{equation}\label{d3}
P(s_n(\bbX)\leq \varepsilon n^{-1/2})\leq \delta \quad\text{for all} \quad n\geq n_0.
\end{equation}
Moreover, since (\ref{d2}) fails to hold for the random sign matrices
($X_{jk}$ being symmetric $\pm 1$ r.v.'s), Spielman and Teng \cite {st}
speculated that for random sign matrices for any $\varepsilon\geq 0$
\begin{equation}\label{d4}
P(s_n(\bbX)\leq \varepsilon n^{-1/2})\leq\varepsilon+c^n\quad 0<c<1.
\end{equation}
Again, (\ref{d4}) has been proved for real random matrices with i.i.d. subgaussian
entries in \cite{rud2}.

 We will adapt Rudelson and Vershynin's
method to obtain the order of the smallest singular
value for complex matrices perturbed by
a constant matrix.

Formally, let $\bbW=\bbX+\bbA_n$, where $\bbA_n$ is a fixed complex matrix
and $\bbX=(X_{jk})$, a random matrix. Denote the singular values
of $\bbW$ by $s_1,\cdots,s_n$ arranged in the non-increasing order.
Particularly, the smallest singular value is
$$
s_n(\bbW)=\inf_{\bbx\in \mathbb{C}^{n}:\| \bbx\|_2=1}\|\bbW
\bbx\|_2,
$$
where $\|\cdot \|_2$ means Euclidean norm,
 and we denote the spectral norm of a matrix by $\|\cdot\|$.

\begin{theorem}
\label{small}Let $\{X_{jk}\}$ be i.i.d. complex r.v.'s
with $EX_{11}=0, \ E|X_{11}|^2=1$ and $E|X_{11}|^3<B$. Let $K\geq
1$.
  Then for every $\e \ge 0$, 
  \begin{equation}              \label{a29}
  P ( s_n(\bbW) \le \e n^{-1/2} )
    \le C\e + c^n+P (\|\bbW\| > K n^{1/2}),
  \end{equation}
where $C> 0$ and $c \in (0,1)$ depend only on $K$, $B$,
$E\big(Re(X_{11})\big)^2$, $E\big(Im(X_{11})\big)^2$,
and $ERe(X_{11})Im(X_{11})$.
\end{theorem}

\begin{remark}
In Theorem \ref{small}, $\e$ is arbitrary. It can depend on $n$.
$K$ is a constant not smaller than $1$. In Section \ref{ccl} when we apply (\ref{a29})
in the proof
of Theorem \ref{theo3}, we will select $\e=n^{-1-\delta}, \ K>4$.
\end{remark}

\begin{remark}
Theorem \ref{small} includes Theorem 5.1 in \cite{rud2} as a special
case, where $\bbA_n=0$, the r.v.'s are real and have
finite fourth moment. Therefore, (\ref{d3}) is true with  $\bbX$ replaced by $\bbW$
when $X_{11}$ has finite fourth moment and $\|\bbA_n\|\leq C\sqrt{n}$ ($C\geq 0$), i.e.,
$$
P(s_n(\bbW)\leq \varepsilon n^{-1/2})\leq \delta \quad\text{for all} \quad n\geq n_0.
$$
Moreover, if $\bbX$ is a subgaussian matrix and $\|\bbA_n\|\leq C\sqrt{n}$,  by Lemma 2.4 of \cite{rud2} or Fact 2.4 of
\cite{rud1}, (\ref{d4}) holds with $\bbX$ replaced by
$\bbW$, i.e.,
$$
P(s_n(\bbW)\leq \varepsilon n^{-1/2})\leq\varepsilon+c^n\quad 0<c<1.
$$
This exponential rate is better than the polynomial rate in Tao and Vu \cite{tao2}.
\end{remark}

Furthermore, for general random matrices, similar to steps (\ref{a34})-(\ref{a35}) in Section 3 one can conclude that
\begin{corollary}\label{cor1}
In addition to the assumptions of Theorem \ref{small}, suppose that $|X_{ij}|\leq \sqrt{n}\varepsilon_n$ and
$\|\bbA_n\|\leq C\sqrt{n}$ with $0\leq C<\infty$, then for any $\varepsilon\geq 0$
\begin{equation}\label{d5}
P ( s_n(\bbW) \le \e n^{-1/2} )
    \le C\e +n^{-l},
\end{equation}
where $l$ is any positive number and $\varepsilon_n\rightarrow0$ with the convergence rate
slower than any preassigned one as $n\to \infty$.
\end{corollary}

\begin{remark}
Taking $\varepsilon=0$, Corollary \ref{cor1} then leads to a polynomial bound for the singularity probability:
$$
P(\bbW_n \ is \ singalur)\leq n^{-l},
$$
with $l$ being any positive number.
\end{remark}
\begin{remark} \label{rem1} For random sign matrices Tao and Vu \cite{tao1} showed that for every $A>0$ there exists $B>0$ so that
$$
P(s_n(\bbX)  \le n^{-B})\leq n^{-A}.
$$
 Recently, Tao and Vu \cite{tao2}
reported a result concerning the smallest singular value of a
perturbed matrix too.  Under some mild conditions, they proved that
$$
P(s_n(\bbW)  \le n^{-B})\leq n^{-A}.
$$
Compared with their results, (\ref{d5}) gives an explicit dependence
between the bound on $s_n(\bbW)$ and probability, while the relationship
between $A$ and $B$ in \cite{tao1} and \cite{tao2} is implicit.
In addition, (\ref{d5}) holds for general random matrices, while
 Tao and Vu's theorem basically applies to discrete random matrices.
\end{remark}

\begin{remark}
In this paper, we will use the letters $B,K_1,K_2$ to denote some finite
absolute constants.
\end{remark}

The argument of Theorem \ref{small} is presented in the next section and
the proof of the circular law is given in the last section.

\section{Smallest singular value}

In this section the smallest singular value of the  matrix
 $\bbX$ perturbed by a constant matrix will be characterized.
 We begin first with the estimation of the so-called small ball
probability.

\subsection{Small ball probability}

The small ball probability is defined as
\begin{equation}\label{a1}
P_\varepsilon(\bbb)=\sup\limits_{v\in \mathbb{C}}
P(|S_n-v|\leq\varepsilon),
\end{equation}
where
\begin{equation}
S_n=\sum\limits_{k=1}^nb_k\eta_k
\end{equation}
with $\eta_1,\cdots,\eta_n$ being i.i.d.
 r.v.'s and $\bbb=(b_1,\cdots,b_n)\in \mathbb{C}^{n}$ (see
\cite{ls}).
If each $\eta_k$ is perturbed by a constant $a_k\in \mathbb{C}$,
then $P_\varepsilon(\bbb)$
does not change, i.e.
\begin{equation}\label{invariance}
P_\varepsilon(\bbb)=\sup\limits_{v\in \mathbb{C}}
P(|\sum\limits_{k=1}^nb_k(\eta_k-a_k)-v|\leq\varepsilon),
\end{equation}

We first establish a small ball probability for big $\varepsilon$
via central limit theorem for complex r.v.'s
$\eta_1,\cdots,\eta_n$. Before we state the next result, let us
introduce some more notation and terminology. $Re(z)$ and $Im(z)$
will denote the real and imaginary part of a complex number $z$.
Write $\eta_{1k}=Re(\eta_k),$ $\eta_{2k}=Im(\eta_k)$,
$\sigma_{1}^2=\sigma_{1k}^2=E(\eta_{1k}-E\eta_{1k})^2, \
 \sigma_{2}^2=\sigma_{2k}^2=E(\eta_{2k}-E\eta_{2k})^2, \
 \sigma_{12}=\sigma_{12k}=E(\eta_{1k}-E\eta_{1k})(\eta_{2k}-E\eta_{2k})$ for
 $k=1,2,\cdots,n$.
For real r.v.'s $\xi$ and $\eta$, if
$\big(E(\xi-E\xi)(\eta-E\eta)\big)^2=E(\xi-E\xi)^2 E(\eta-E\eta)^2>0$,
then we will say
that $\xi$ and $\eta$ are linearly correlated.

\begin{theorem}\label{theo1}
Let $\eta_1,\cdots,\eta_n$ be i.i.d. complex r.v.'s
with variances at least $1$, $E|\eta_1|^3<B$ and let $b_1,\cdots,b_n$ be
complex numbers such that $0<K_1\leq|b_k|\leq K_2$ for all $k$. Then
for every $\varepsilon>0$,
\begin{equation}\label{a50}
P_\varepsilon(\bbb)\leq
\frac{C}{\sqrt{n}}\left(\frac{\varepsilon}{K_1}+(\frac{K_2}{K_1})^3\right),
\end{equation}
where $C$ is a finite constant depending only on
 $B$, $\sigma_{1},\ \sigma_{2}$ and $\sigma_{12}$.
\end{theorem}

\begin{proof}
Suppose first that $Re(\eta_k)$ and $Im(\eta_k)$ are linearly correlated,
$k=1,\cdots,n$.
Then $\eta_k-E\eta_k=\xi_k(1+ib_0)/(1+b_0^2)^{1/2}$ $a.s.$, where
$\xi_k=(1+b_0^2)^{1/2}Re(\eta_k-E\eta_k)$ and
$b_0$ is an absolute real constant. Write
$\tilde b_k=b_k(1+ib_0)/(1+b_0^2)^{1/2}$ which
satisfies $K_1\leq|\tilde b_k|\leq K_2$. Let $\tilde b_{1k}=Re(\tilde b_k)$ and
$\tilde b_{2k}=Im(\tilde b_k)$. Noting that
$$\sup\limits_{v\in \mathbb{C}}P(|S_n-v|\leq \varepsilon) \leq
\sup\limits_{v\in \mathbb{C}}
P(|\sum_{k=1}^n\tilde b_{1k}\xi_k-Re(v)|\leq \varepsilon, \
|\sum_{k=1}^n\tilde b_{2k}\xi_k-Im(v)|\leq \varepsilon)
$$
and either $\sum_{k=1}^n\tilde b_{1k}^2\geq nK_1^2/2$ or
$\sum_{k=1}^n\tilde b_{2k}^2\geq nK_1^2/2$,  we can complete the
proof for the linearly correlated
case by Berry-Esseen inequality.

The case where $Re(\eta_k)=0$ or $Im(\eta_k)=0$ $a.s.$ follows from
Berry-Esseen inequality directly.

Now suppose $Re(\eta_k)$ and $Im(\eta_k)$ are not linearly correlated,
and $P\big(Re(\eta_k)=0\big)<1$, $P\big(Im(\eta_k)=0\big)<1$.
 Let $b_k=b_{1k}+ib_{2k}$ and
$v=v_1+iv_2$. Define
$\hat{\eta}_{1k}=b_{1k}\eta_{1k}-b_{2k}\eta_{2k}$ and
$\hat{\eta}_{2k}=b_{1k}\eta_{2k}+b_{2k}\eta_{1k}$. Obviously,
$\sum\limits_{k=1}^nE|\hat{\eta}_{jk}-E\hat{\eta}_{jk}|^3\leq
\sum\limits_{k=1}^n E|b_k(\eta_k-E\eta_k)|^3 \leq
8B\|b\|_3^3,j=1,2$, where $\|b\|_3^3=\sum\limits_{k=1}^n|b_k|^3$.
In order to apply Berry-Esseen inequality, we need to get a lower bound for
$E|\hat{\eta}_{jk}-E\hat{\eta}_{jk}|^2$. For $j=1$, we have
\begin{eqnarray*}
&&E|\hat{\eta}_{1k}-E\hat{\eta}_{1k}|^2 \\
&=&b_{1k}^2\sigma_{1k}^2+b_{2k}^2\sigma_{2k}^2-2b_{1k}b_{2k}\sigma_{12k} \\
&=&|b_k|^2\Big((|b_{1k}|\sigma_{1k}/|b_k|-|b_{2k}|\sigma_{2k}/|b_k|)^2 \\
&& \qquad +2|b_{1k}b_{2k}||b_k|^{-2}\big(\sigma_{1k}\sigma_{2k}-\sign(b_{1k}b_{2k})
\sigma_{12k}\big)\Big).
\end{eqnarray*}
For $t\in [0,1]$, let $f(t)=(t\sigma_{1}-\sqrt{1-t^2}\sigma_{2})^2+2t\sqrt{1-t^2}
(\sigma_{1}\sigma_{2}\pm\sigma_{12})$. So the smallest value
$a=\min_{t\in[0,1]}f(t)$ of $f(t)$ in
$[0,1]$ is attainted at $0$ or $1$ or some $t_0\in (0,1)$.
Therefore, $a$ is a positive
constant depending only on $\sigma_{1},\ \sigma_{2}$ and $\sigma_{12}$.
Hence $E|\hat{\eta}_{1k}-E\hat{\eta}_{1k}|^2 \geq a|b_k|^2$.
Similarly, $E|\hat{\eta}_{2k}-E\hat{\eta}_{2k}|^2 \geq a|b_k|^2$.
By Berry-Esseen inequality, one can then conclude
that
\begin{equation}
\label{a3} \sup\limits_{v_1\in\mathbb{R}}
P(|\sum\limits_{k=1}^n(\hat{\eta}_{1k}-E\hat{\eta}_{1k})-v_1|\leq
\frac{\varepsilon}{\sqrt{2}})\leq
\frac{C\varepsilon}{\|\bbb\|_2}+C\left(\frac{\|\bbb\|_3}{\|\bbb\|_2}\right)^3
\end{equation}
and
\begin{equation}
\label{a4}
\sup\limits_{v_2\in\mathbb{R}}P(|\sum\limits_{k=1}^n(\hat{\eta}_{2k}
-E\hat{\eta}_{2k})-v_2|\leq\frac{\varepsilon}{\sqrt{2}})\leq
\frac{C\varepsilon}{\|\bbb\|_2}+C\left(\frac{\|\bbb\|_3}{\|\bbb\|_2}\right)^3,
\end{equation}
where $C$ is a constant depending only on $B$,
$\sigma_{1},\ \sigma_{2}$ and $\sigma_{12}$.

 Thus (\ref{a50}) follows from (\ref{a3}), (\ref{a4}) and the following inequality
\begin{eqnarray*}
\sup\limits_{v\in \mathbb{C}}
P(|S_n-v|\leq\varepsilon)
&\leq&\sup\limits_{v_1\in
\mathbb{R}}P(|\sum\limits_{k=1}^n(\hat{\eta}_{1k}-E\hat{\eta}_{1k})
-v_1|\leq\frac{\varepsilon}{\sqrt{2}}) \\
& &  \qquad +
\sup\limits_{v_2\in
\mathbb{R}}P(|\sum\limits_{k=1}^n(\hat{\eta}_{2k}-E\hat{\eta}_{2k})
-v_2|\leq\frac{\varepsilon}{\sqrt{2}}).
\end{eqnarray*}
\end{proof}

Theorem \ref{theo1} only yields a polynomial rate $n^{-1/2}$.
Next, an improved small ball probability is needed for our future
use. To this end, some concepts will be presented which are parallel
to those of \cite{rud2}.

Denote the unit sphere in $\mathbb{C}^{n}$ by $S^{n-1}$.

\begin{definition}
Let $\alpha\in(0,1)$ and $\tau\geq 0$. The essential least common
denominator of a vector $\bbb\in\mathbb{C}^n$, denoted by
$D(\bbb)=D_{\alpha,\tau}(\bbb)$, is defined to be the infimum of
$t>0$ so that all coordinates of the vector $t\bbb$ are of distance
at most $\alpha$ from nonzero integers except $\tau$ coordinates.
\end{definition}

\begin{definition}
Suppose that $\gamma,\rho\in (0,1)$. A vector $\bbb\in\mathbb{C}^n$
is sparse if $|supp(\bbb)|\leq \gamma n$. A vector $\bbb\in S^{n-1}$
is compressible if $\bbb$ is within Euclidean distance $\rho$ from the
set of all sparse vectors. All vectors $\bbb\in S^{n-1}$ except
compressible vectors are called incompressible. Let
$Sparse=Sparse(\gamma), Comp=Comp(\gamma,\rho)$ and
$Incomp=Incomp(\gamma,\rho)$ denote, respectively, the sets of
sparse, compressible and incompressible vectors.
\end{definition}

\begin{definition}
For some $K_1,K_2>0$, the spread part of a vector
$\bbb\in\mathbb{C}^n$ is defined as
$$
\hat{\bbb}=(\sqrt{n}b_k)_{k\in\sigma(\bbb)},
$$
where the subset $\sigma(\bbb)\subseteq\{1,\cdots,n\}$ is given by
$\{k:\quad K_1\leq\sqrt{n}|b_k|\leq K_2\}$. Similarly, for $j=1,2$,
define
$$
\hat{\bbb}_j=(\sqrt{n}b_{jk})_{k\in\sigma(\bbb)},\quad
|\hat{\bbb}_j|=(\sqrt{n}|b_{jk}|)_{k\in\sigma(\bbb)}, \quad
\hat{|\bbb|}=(\sqrt{n}|b_{k}|)_{k\in\sigma(\bbb)},
$$
where $b_{1k}$ and $b_{2k}$ denote, respectively, the real part and
imaginary part of $b_k$.
\end{definition}

Similar to the real case, the complex incompressible vector are also
evenly spread, i.e. many coordinates are of the order $n^{-1/2}$.
\begin{lemma}
\label{lem1} Let $\bbb\in Incomp(\gamma,\rho)$. Then there is a
set $\sigma_1(\bbb)\subset\{1,\cdots,n\}$ of cardinality
$|\sigma_1(\bbb)|\geq c n$ with $c\geq \rho^2\gamma /4$ so that for
$j=1$ or $2$,
\begin{equation}\label{a8}
\frac{\rho}{2\sqrt{2n}}\leq |b_{jk}|\leq\frac{1}{\sqrt{\gamma n}}
\ \mbox{for all}\  k\in \sigma_1(\bbb).
\end{equation}
\end{lemma}

\begin{proof} By Lemma 3.4 in \cite {rud2}, for $\bbb\in Incomp(\gamma,\rho)$,
there is a set $\sigma(\bbb)$ of cardinality
$|\sigma(\bbb)|\geq\frac{1}{2}\rho^2\gamma n$ so that
\begin{equation*}
\frac{\rho}{\sqrt{2n}}\leq |b_k|\leq\frac{1}{\sqrt{\gamma
n}}\quad\mbox{for all $k\in\sigma(\bbb)$}.
\end{equation*}
Hence $|b_{1k}|\leq 1/\sqrt{\gamma n}$ and $|b_{2k}|\leq 1/\sqrt{\gamma n}$ if
$k\in\sigma(\bbb)$. On the other hand, either $b_{1k}$ or $b_{2k}$ must be bigger than
$\rho(2\sqrt{2n})^{-1}$. The assertion follows.
\end{proof}

The following result refines Theorem \ref{theo1}.
\begin{theorem}
\label{theo2} Let $\bbb=(b_1,\cdots,b_n)\in \mathbb{C}^n$ whose spread part
$\hat \bbb$ is well defined (for some fixed truncation levels $K_1, \ K_2>0$). Suppose
$0<\alpha<K_1/6K_2$ and $0<\beta<1/2$.

 (1) Suppose that
$\eta_1,\cdots,\eta_n$ are i.i.d. real r.v.'s, or
imaginary r.v.'s, or complex ones with linearly correlated $Re (\eta_k)$
and $ Im (\eta_k), k=1,2,\cdots,n$. If $E|\eta_k-E\eta_k|^2=1$
and $E|\eta_k|^3<B$, for any $\varepsilon\geq0$, then
\begin{equation}\label{a5}
P_\varepsilon(\bbb)\leq
\frac{C}{\sqrt{\beta}}\left(\varepsilon+\frac{1}{\sqrt{n}
\max \{D_{\alpha,\beta n}(\hat{\bbb}_1),
D_{\alpha,\beta n}(\hat{\bbb}_2)\}}\right)+C\exp(-c\alpha^2\beta n),
\end{equation}
where $C,c>0$ depend only on $B,K_1,K_2$.

 (2) Let $\eta_1,\cdots,\eta_n$ be
i.i.d. complex r.v.'s with $E|\eta_k-E\eta_k|^2=1$ and
$E|\eta_k|^3<B$, then (\ref{a5}) holds or
\begin{equation}\label{a16}
P_\varepsilon(\bbb)\leq
\frac{C}{\sqrt{\beta}}\left(\varepsilon+\frac{1}{\sqrt{n}D_{\alpha,\beta
n}(\hat{|\bbb|})}\right)+C\exp(-c\alpha^2\beta n)
\end{equation}
where $C,c>0$ depend only on $B, \ K_1, \ K_2$,
 $\sigma_{1},\ \sigma_{2}$ and $\sigma_{12}$.
\end{theorem}


\begin{proof}
Since $P_\varepsilon(\bbb)=\sup\limits_{v\in \mathbb{C}}
P(|S_n-ES_n-v|\leq\varepsilon)$, we can assume that $E\eta_k=0$.

(1). We only consider the case where the r.v.'s $\{\eta_k\}$ are
real. The other two cases follow from the real case.
 Let $b_k=b_{1k}+ib_{2k}$ and $v=v_1+iv_2$.
Noting that
\begin{eqnarray*}
&&\sup\limits_{v\in \mathbb{C}}P(|S_n-v|\leq \varepsilon)\\
&\leq& \min\Big(\sup\limits_{v_1\in \mathbb{R}}
P(|\sum_{k=1}^n b_{1k}\eta_k-v_1|\leq \varepsilon),
\sup\limits_{v_2\in \mathbb{R}}P(|\sum_{k=1}^n  b_{2k}\eta_k
-v_2|\leq \varepsilon)\Big)
\end{eqnarray*}
Then Corollary 4.9 in \cite{rud2} leads to (\ref{a5}).

(2). 
For the moment we assume that
$$
1\leq |b_k|\leq K \quad\mbox{for all k }.
$$
Let $b_k=b_{1k}+ib_{2k}$, $\eta_k=\eta_{1k}+i\eta_{2k}$ and
$v=v_1+iv_2$. It is observed that Theorem \ref{theo1} implies
Theorem \ref{theo2} for big values of $\varepsilon$ (constant order
or even larger). Therefore we can suppose in what follows that
$$
\varepsilon\leq l_1,
$$
where $l_1$ is a constant which will be specified later.

If the real part of $\eta_1$ is linearly correlated to the imaginary part of
$\eta_1$, then we have (\ref{a5}). 
 Therefore we
assume in the sequel that $\eta_{11}$ is not linearly correlated to
$\eta_{21}$. 

Set $\zeta_k=\frac{1}{|b_k|}|\xi_k-\xi_k'|$ where
$\xi_k=b_{1k}\eta_{1k}-b_{2k}\eta_{2k}$ and $\xi_k'$ is an
independent copy of $\xi_k$. Then
\begin{eqnarray}
\frac{1}{2}E\zeta_k^2&=&\frac{1}{|b_k|^2}E|b_{1k}\eta_{1k}-b_{2k}\eta_{2k}|^2.
\label{a11}
\end{eqnarray}
As in the proof of Theorem \ref{theo1}
$$
E\zeta_k^2\geq 2a>0,
$$
where $a$ is some positive
constant depending only on $\sigma_{1},\ \sigma_{2}$ and $\sigma_{12}$.

On the other hand,
$E\zeta_k^3\leq 64B$. The Paley-Zygmund inequality (\cite{rud1})
gives that
$$
P(\zeta_k>\sqrt{a})\geq
\frac{(E\zeta_k^2-a)^3}{(E\zeta_k^3)^2}\geq\frac{a^3}{64^2B^2}=:\beta,
$$
which is a positive constant depending only on $B$, $\sigma_1$, $\sigma_2$
and $\sigma_{12}$. Following \cite{rud2} we introduce a new
r.v. $\hat{\zeta}_k$ conditioned on
$\zeta_k>\sqrt{a}$, that is, for any measurable function $g$
$$
Eg(\hat{\zeta}_k)=\frac{Eg(\zeta_k)I(\zeta_k>\sqrt{a})}{P(\zeta_k>\sqrt{a})},
$$
which entails
 \begin{equation}
\label{a13} Eg(\zeta_k)\geq\beta
Eg(\hat{\zeta}_k).
\end{equation}

From Esseen inequality, one has
\begin{eqnarray}
P_\varepsilon(\bbb)&\leq&
\sup\limits_{v_1\in \mathbb{R}}P(|\sum\limits_{k=1}^n\xi_k-v_1|\leq\varepsilon)\non
&\leq& C\int^{\pi/2}_{-\pi/2}|\phi(t/\varepsilon)|dt,\label{a14}
\end{eqnarray}
where
$$
\phi(t):=E\exp(i\sum\limits_{k=1}^n\xi_kt).
$$
With the notation $\phi_k(t)=E\exp(i\xi_kt)$, it is observed that
$$|\phi_k(t)|^2=E\cos(|b_k|\zeta_kt),$$ and we then have
\begin{align*}
|\phi(t)|
  &\le \prod_{k=1}^n \exp \Big( - \frac{1}{2} (1 - |\phi_k(t)|^2) \Big)\\
  &= \exp \Big( -E \sum_{k=1}^n \frac{1}{2} (1 - \cos (|b_k| \zeta_k t)) \Big)
  = \exp \big( -E g(\z_k t) \big),
\end{align*}
where
$$
g(t) := \sum_{k=1}^n \sin^2 \big( \frac{1}{2} |b_k|  t \big).
$$
This, together with (\ref{a13}), gives
$$
|\phi(t)| \le \exp \big( - \b \; E g(\hat{\z}_k t) \big).
$$
Consequently, (\ref{a14}) becomes
\begin{align}                           \label{SBP via int}
P_\e(\bbb)
  &\le C \int_{-\pi/2}^{\pi/2} \exp \big( - \b \; E g(\hat{\z}_k t / \e) \big) \; dt \nonumber\\
  &\le C E \; \int_{-\pi/2}^{\pi/2} \exp \big( - \b g(\hat{\z}_k t / \e) \big) \; dt \nonumber\\
  &\le C \sup_{z \ge \sqrt{a}} \;
    \int_{-\pi/2}^{\pi/2} \exp \big( - \b g(zt/\e) \big) \; dt.
\end{align}

Let
$$ M := \max_{|t| \le \pi/2} g(zt/\e) = \max_{|t| \le \pi/2}
\sum_{k=1}^n \sin^2 (|b_k| z t / 2 \e)
$$
and the level sets of $g$ be
$$
T(m,r) := \{ t: \; |t| \le r, \; g(z t/\e) \le m \}.
$$
As in \cite{rud2}, one can prove that
 $$
  \frac{n}{4}  \le M \le n,
  $$
  by taking $\varepsilon<(\pi\sqrt{a})/4=l_1$. All the remaining arguments
   including the analysis for the
  level sets $T(m,r)$ are similar to those of \cite{rud2} and so we here omit
  the details. Thus, one can conclude that
 for every $\varepsilon\geq0$
\begin{equation}\label{a15}
|P_\varepsilon(\bbb)|\leq
\frac{C}{\sqrt{\tau}}\left(\varepsilon+\frac{1}{D_{\alpha,\tau}(|\bbb|)}\right)+C\exp(-\frac{c\alpha^2\tau}{A^2}).
\end{equation}
where $0<\tau<n$, $|\bbb|=(|b_1|,\cdots,|b_n|)$ and $C,c>0$ are positive
constants depending only on $B$, $\sigma_1$, $\sigma_2$
and $\sigma_{12}$..

Finally, combining (\ref{a15}) and Lemma 2.1 in \cite{rud2} one can
obtain the small ball probability for complex case (when applying
(\ref{a15}) to the spread part of the vector $\bbb$ one can suppose
that $K_1=1$ by re-scaling $b_k$ and $\alpha$). Thus we complete the
proof.
\end{proof}

To treat the compressible vector, the following lemma is needed.

\begin{lemma}\label{lem2}
Suppose that $\eta_1,\cdots,\eta_n$ are i.i.d. centered complex
r.v.'s with $E|\eta_k|^2=1$ and $E|\eta_k|^3\leq B$. Let
$\{a_{jk},j,k=1,\cdots,n\}$ be complex numbers.
 Then
for $0<\lambda<1$ and any vector $\bbb=(b_1,\cdots,b_n)\in S^{n-1}$ there is
$\mu\in(0,1)$ such that the sum
$S_{nj}=\sum\limits_{k=1}^nb_k(\eta_k-a_{jk})$ satisfy
$$
P(|S_{nj}|>\lambda)\geq\mu
$$
where $\mu$ depends only on $\lambda$ and $B$.
\end{lemma}

\begin{proof}  Simple calculation indicates that
$$
E|S_{nj}|^2=|\sum\limits_{k=1}^nb_ka_{jk}|^2+1.
$$
On the other hand by Burkholder inequality we have
\begin{eqnarray*}
E|S_{n1}|^3 &\leq&
4\Big(|\sum\limits_{k=1}^nb_ka_{jk}|^3+E|\sum\limits_{k=1}^nb_k\eta_k|^3\Big)\\
&\leq& C\Big(|\sum\limits_{k=1}^nb_ka_{jk}|^3+(\sum\limits_{k=1}^n|b_k|^2E|\eta_k|^2)^{3/2}+\sum\limits_{k=1}^n|b_k|^3E|\eta_k|^3\Big)\\
&\leq& C\left(|\sum\limits_{k=1}^nb_ka_{jk}|^3+1+B\right).
\end{eqnarray*}
Hence Paley-Zygmund inequality gives that
$$
P(|S_{nj}|>\lambda)\geq\frac{(E|S_{nj}|^2-\lambda^2)^3}{(ES_{nj}^3)^{2}}\geq
\frac{(c_{nj}^2+1-\lambda^2)^3}{C(c_{nj}^3+1+B)^2},
$$
where
$$
c_{nj}=|\sum\limits_{k=1}^nb_ka_{jk}|.
$$
Take
$$
f(t)=\frac{(t^2+1-\lambda^2)^3}{(t^3+1+B)^2}, \quad
t\in(0,\infty).
$$
Then one can conclude that
$$
\mu:=\min\limits_{t\in(0,\infty)}f(t)>0
$$
and then
$$
P(|S_{nj}|>\lambda)\geq \mu>0
$$
where $\mu$ depends only on $\lambda$ and $B$.
\end{proof}

\subsection{Proof of Theorem \ref{small}}

The whole argument is similar to that of \cite {rud2} and we only
sketch the proof. For more details one can refer
to \cite{rud2}.

Since $ S^{n-1}$ can be decomposed as the union of
$Comp$ and $Incomp$, we then consider the smallest singular value on
 each set separately.

By Lemma \ref{lem2} there are $c_1>0$ and $v\in(0,1)$ depending on $\mu$
only so that
$$
P(\|\bbW \bbb\|_2<c_1\sqrt{n})\leq v^n, \ \bbb \in S^{n-1}.
$$
Actually, the proof is similar to that of Proposition 3.4 in
\cite{rud1}. The only difference is that we should use our
Lemma \ref{lem2} instead of Lemma 3.6 in \cite{rud1}.
 Therefore similar to Lemma 3.3 in \cite{rud2}, we have,
there exist $\gamma,\rho,c_2,c_3>0$ so that
\begin{equation}\label{a17}
  P \big( \inf_{\bbb \in \Comp(\gamma,\rho)} \|\bbW\bbb\|_2 \le c_2 n^{1/2}
   \big)
  \le e^{-c_3 n}+P(\|\bbW\| > K n^{1/2}),
  \end{equation}
where $K\geq 1$.

Let $\bbX_1,\cdots,\bbX_n$ denote the column vectors of $\bbW$ and
$H_k$ the span of all columns except the $k$-th column. One can
check that Lemma 3.5 in \cite{rud2} is still true in complex case
and hence
\begin{eqnarray}
 P(\inf\limits_{\bbb\in Incomp(\gamma,\rho)}\|\bbW \bbb\|_2\leq\varepsilon\rho n^{-1/2})
& \leq&\frac{1}{\gamma n}\sum_{k=1}^nP ( \dist(\bbX_k, H_k) <
\e)\non &\leq&\frac{1}{\gamma n}\sum_{k=1}^nP (|\langle \bbY_k,
\bbX_k\rangle|< \e),\label{a18}
\end{eqnarray}
where $\bbY_k$ is any unit vector orthogonal to $H_k$ and can be
chosen to be independent of $\bbX_k$. Here
$\langle\cdot,\cdot\rangle$ is the canonical inner product in
$\mathbb{C}^n$.

  When all $\{X_{jk}\}$ are real r.v.'s,
or when $Re(X_{jk})$ and $Im(X_{jk})$ are linearly correlated
 or when $Re(X_{jk})=0$ we have
\begin{eqnarray}
P (|\langle \bbY_k,\bbX_k\rangle|< \e \text{ and } U_{K})
 &\leq & P
\big( \bbY_{k}\in Comp \text{ and } U_K\big) \non
 &\ &+P(|\langle \bbY_k,\bbX_k\rangle|< \e, \bbY_{k}\in Incomp \text{ and } U_K),
      \non  \label{a19}
\end{eqnarray}
where $U_K$ denotes the event that $\|\bbW\| \le K n^{1/2}$. One can
check that Lemma 3.6 in \cite{rud2} applies to complex case and
hence
$$
P(\bbY_{k}\in Comp\ \text{and} \ U_K)\leq e^{-c_4n},
$$
where $c_4$ is a constant depending only on $B, \  K$,
 $\sigma_{1},\ \sigma_{2}$ and $\sigma_{12}$. Further,
\begin{eqnarray*}
&&P(|\langle \bbY_k,\bbX_k\rangle|< \e, \ \bbY_{k}\in Incomp \text{ and } U_K)\qquad \\
&\leq& \sum_{j=1}^2P \left(V_{jk}, \ U_K,
      \ D_{\a, \b n}(\hat{\bbY}_{jk}) < e^{cn}
       \text{ and } \bbY_{k}\in Incomp \right)\\
&+&\sum_{j=1}^2E\left[I\big(D_{\a, \b n}(\hat{\bbY}_{jk}) \geq e^{cn}\text{ and } \bbY_{k}\in Incomp\big)
P \left(|\< \bbY_k,\bbX_k\> | < \e
        |\bbY_k\right)\right]\qquad \\
\end{eqnarray*}
where $V_{1k}$ and $V_{2k}$ denote, respectively, the events that
 the real part and imaginary part of the vector
$\bbY_{k}\in Incomp$ satisfy (\ref{a8}) in Lemma \ref{lem1}, $\hat{\bbY}_{1k}$ and
$\hat{\bbY}_{2k}$ denote, respectively, the spread part of the real
part and imaginary part of the vector $\bbY_{k}$. By (\ref{a5}) in Theorem
\ref{theo2} and (\ref{invariance}) we have
$$
I\big(D_{\a, \b n}(\hat{\bbY}_{jk}) \ge e^{cn}\big)P \left( |\< \bbY_k,\bbX_k\> | < \e
       |\bbY_k\right)
      \leq c_5\varepsilon+c_6e^{-c_7n},
$$
where $c_5,c_6,c_7$ are positive constants depending only on $B$,
 $\sigma_{1},\ \sigma_{2}$ and $\sigma_{12}$.

On the other hand,
\begin{eqnarray*}
&&P \left(V_{1k},\  U_K,
      \  D_{\a, \b n}(\hat{\bbY}_{1k}) < e^{cn}
       \text{ and } \bbY_{k}\in Incomp \right)\\
      &\leq& \sum_{D \in \DD} \P(\bbY_k \in S_D, \  U_K \text{ and } V_{1k}).
\end{eqnarray*}
Here the level set $S_D \subseteq S^{n-1}$ is defined as
  $$
  S_D := \{ \bbY_k \in \Incomp : \; D \le D_{\a, n_0/2}(\hat{\bbY}_{1k}) < 2D \}.
  $$
  and
$$
\DD = \{ D :\; D_0 \le D < e^{cn}, \; D = 2^k, \; k \in \Z \},
$$
where $\alpha$ and $D_0$ are some constants. For more details about $\alpha$ and $D_0$,
see \cite{rud2}. Further, one can similarly prove
that Lemma 5.8 in \cite{rud2} holds in our case and therefore we
obtain
$$
\P(\bbY_k \in S_D \text{ and } U_K )\leq e^{-n},
$$
which, combined with the fact that the cardinal number $|\DD|$ is of order $n$, then implies that
$$
P \left(V_{1k},\  U_K,
      \  D_{\a, \b n}(\hat{\bbY}_{1k}) < e^{cn}
       \text{ and } \bbY_{k}\in Incomp \right)\leq e^{-c_8n},
$$
where $c_8>0$. Similarly, one may also show that
$$
P \left(V_{2k}, \ U_K,
      \  D_{\a, \b n}(\hat{\bbY}_{2k}) < e^{cn}
       \text{ and } \bbY_{k}\in Incomp \right)\leq e^{-c_8n}.
$$
Picking up the above argument one can conclude that
$$
P \big( |\< \bbY_k, \bbX_k\> | < \e \text{ and } \|\bbW\| \le K
n^{1/2} \big) \le C \e + e^{-c' n},
$$
which further gives that
\begin{equation}
  P \big( \inf_{\bbx \in \Incomp(\gamma,\rho)} \|\bbW \bbx\|_2 \le \e \rho n^{-1/2} \big)
  \le \frac{C}{\d} (\e + c^n) + P (\|\bbW\| > K n^{1/2}),\label{a20}
  \end{equation}
where $C>0$ and  $c\in(0,1)$ depend only on $K$, $B$,
 $\sigma_{1},\ \sigma_{2}$ and $\sigma_{12}$.

For all the remaining case, i.e. $Re(X_{jk})Im(X_{jk})\not\equiv0$, and
$Re(X_{jk})$, $Im(X_{jk})$ are not linearly correlated, one has
\begin{eqnarray}
&&P (|\langle \bbY_k,\bbX_k\rangle|< \e \text{ and } U_{k})
 \leq  P
\big( D_{\a, \b n}(|\bbY_k|) < e^{cn}
      \text{ and } U_K\big) \non
 &&\quad \qquad + E\left[I\big(D_{\a, \b n}(|\bbY_k|) \ge e^{cn}\big)P \left( |\< \bbY_k,\bbX_k\> | < \e
      \Big|\bbY_k\right)
      \right],
      \non
 \label{a26}
\end{eqnarray}
and one can similarly obtain (\ref{a20}) for complex case. Theorem \ref{small}
follows from (\ref{a17})-(\ref{a26}) immediately.

\section{The convergence of logarithmic potential and circular law} \label{ccl}

In this part the logarithmic potential will be used to show that the
circular law is true. According
to Lower Envelop Theorem and Unicity Theorem (see Theorem 6.9, p.73,
and Corollary 2.2, p.98, in \cite{saff}), it suffices to show that the corresponding
potential converges to the potential of the circular law.

To make use of Theorem \ref{small} one needs to bound the maximum
singular value of $\bbW$. To this end, we would like to present an
important fact which was proved in \cite{yin1}, that is, if (1)
$EX_{jk}=0$, (2) $|X_{jk}|\leq\sqrt{n}\varepsilon_n,$ (3) $
E|X_{jk}|^2\leq 1 \ \text{and}\  1\geq E|X_{jk}|^2\rightarrow1 $ and
(4) $ E|X_{jk}|^l\leq c (\sqrt{n}\varepsilon_n)^{l-3}\  \text{for}\
l\geq 3$, where $\varepsilon_n\rightarrow0$ with the convergence rate
slower than any preassigned one as $n\to \infty$. Then
for any $K>4$
\begin{equation}\label{a21}
P (\|\bbX\bbX^*\| > K n)=o(n^{-l}),
\end{equation}
where $l$ is any positive number (proved  for real
case in \cite{yin1}, for complex case see Chapter 5 of \cite{bai3}).

Let the random matrix $\hat{\bbX}=(\hat{X}_{jk})$ with
$\hat{X}_{jk}=X_{jk}I(|X_{jk}|\leq \sqrt{n}\varepsilon_n)$. Then one
can show that
\begin{equation}\label{a33}
P(\hat{\bbX}\neq \bbX, i.o.)=0,
\end{equation}
see Lemma 2.2 of \cite{yin1}  (the argument of the complex case is
similar to that of the real one). Here the notation $i.o.$ means
infinitely often. Thus it is sufficient to consider the random
matrix $\hat{\bbX}$ in order to prove the conjecture.

Taking $\bbA_n=E\hat{\bbX}-z\sqrt{n}\bbI$ in Theorem \ref{small}
one can obtain that
\begin{equation}              \label{a34}
  P ( s_n(\hat{\bbX}-z\sqrt{n}\bbI) \le \e n^{-1/2} )
    \le C\e + c^n+P (\|\hat{\bbX}-z\sqrt{n}\bbI\| > K n^{1/2}),
  \end{equation}
  where $E\hat{\bbX}=(E\hat X_{kj})$. Here one should note that from (\ref{a34})
  re-scaling  the underlying r.v.'s is trivial. Moreover
  $$
\|E\hat{\bbX}-z\sqrt{n}\bbI\|\leq
|z|\sqrt{n}+\frac{E|X_{11}|^4}{n\varepsilon_n^3}.
  $$
Therefore, applying (\ref{a21}) and choosing an appropriate $K$ in
(\ref{a34}), we have
\begin{equation}    \label{a35}
  P ( s_n(\hat{\bbX}-z\sqrt{n}\bbI) \le \e n^{-1/2} )
    \le C\e + c^n+n^{-l}
    \end{equation}
where both $C> 0$ and $c \in (0,1)$ depend only on $K$, $E|X_{11}|^3$,
 $E\big(Re(X_{11})\big)^2$, $E\big(Im(X_{11})\big)^2$,
and $ERe(X_{11})Im(X_{11})$.

In the sequel, to simplify the notation, we still use the notation
$\bbX$ instead of $\hat{\bbX}$ and $\mu_n(x,y)$ instead of the
empirical spectral distribution corresponding to $\hat{\bbX}$. But
one should keep in mind that $\{X_{kj}\}$ are non-centered and
$|X_{kj}|\leq \sqrt{n}\varepsilon_n$.

Let
$$\bbH_n=(n^{-1/2}\bbX-z\bbI)(n^{-1/2}\bbX-z\bbI)^*$$
for each $z=s+it\in\mathbb{C}$. Here $(\cdot)^*$ denotes the
transpose and complex conjugate of a matrix. Let $v_n(x,z)$ be the
empirical spectral distribution of Hermitian matrix $\bbH_n$.

Before we prove the convergence of the logarithmic potential of $\mu_n(x,y)$,
we will characterize the relation
between the potential of the circular law $\mu(x,y)$ and the
integral of logarithmic function with respect to $v(x,z)$, the
limiting distribution of $v_n(x,z)$  as below.

\begin{lemma}\label{lem4}
$$
\int\int\log\frac{1}{|x+iy-z|}d\mu(x,y)=-\frac{1}{2}\int^\infty_0\log
x v(dx,z).
$$
\end{lemma}

\begin{proof} Let $x+iy=re^{i\theta},r>0$. One can then verify that
\begin{equation}\label{a27}
\int^{\pi}_{-\pi}\log|z-re^{i\theta}|d\theta=
 \begin{cases}2\pi\log r& \text{if $|z|\leq r$},\\
              2 \pi\log |z| & \text{if $|z|>r$}.
\end{cases}
\end{equation}
It follows that
\begin{equation}\label{a28}
\int\int\log\frac{1}{|x+iy-z|}d\mu(x,y)=
\begin{cases}
 2^{-1}(1-|z|^2) & \text{if $|z|\leq 1$},\\
                      -\log|z| & \text{if $|z|>1$}.
\end{cases}
\end{equation}

On the other hand by Lemma 4.4 in \cite{bai2} one has
$$
\frac{d}{ds}\int^\infty_0\log x v(dx,z)=g(s,t),
$$
where
\begin{eqnarray*}
 g(s,t)=
 \begin{cases}
 \frac{2s}{s^2+t^2}& \text{if $s^2+t^2>1$}\\
              2s   & \text{otherwise}.
  \end{cases}
\end{eqnarray*}
Therefore for any $z=s+it,z_1=s_1+it$ with $|z_1|>1$, we have
\begin{equation}
\int^\infty_0\log x v(dx,z)-\int^\infty_0\log x
v(dx,z_1)+\log|z_1|^2=\int^s_{s_1}g(u,t)du+\log|z_1|^2.
\end{equation}
Let $s_1\rightarrow\infty$ and then $|z_1|\rightarrow\infty$.
Therefore, from Lemma 4.2 of \cite{bai2} the left and right end
point, $x_1$ and $x_2$, of the support of $v(\cdot,z_1)$ satisfy
$$
\frac{x_j}{|z_1|^2}=1+o(1),\quad j=1,2,
$$
which implies that
$$
\int^\infty_0\log x
v(dx,z_1)-\log|z_1|^2=\int^{x_2}_{x_1}\log\frac{x}{|z_1|^2}v(dx,z_1)\rightarrow
0,
$$
as $s_1\rightarrow\infty$. In addition,
\begin{equation}
\int^s_{s_1}g(u,t)du+\log|z_1|^2=
\begin{cases}
|z|^2-1 &  \text{if $|z|\leq1$} \\
 \log|z|^2 & \text{if $|z|>1$}.
\end{cases}
\end{equation}
Thus Lemma \ref{lem4} is complete.
\end{proof}

We now proceed to prove the convergence of the potential of
$\mu_n(x,y)$. The potential of $\mu_n(x,y)$ is
\begin{eqnarray}
U^{\mu_n(x,y)}
&=&-\frac{1}{n}\log\Big|\det\left(n^{-1/2}\bbX-z\bbI\right)\Big|\non
&=&-\frac{1}{2n}\log\Big|\det(\bbH_n)\Big|\non
&=&-\frac{1}{2}\int^\infty_0\log x v_n(dx,z),
\end{eqnarray}
where $\bbI$ is the identity matrix.
We will prove
$$\int^\infty_0\log x v_n(dx,z) \stackrel{a.s.}\longrightarrow
 \int^\infty_0\log x v(dx,z)$$ as
$n\to \infty$. Observe that by the fourth moment condition
$$
\lambda_{\max}(\bbH_n)\leq 2(\lambda_{\max}(n^{-1}\bbX\bbX^*)+|z|^2)\stackrel{a.s.}\longrightarrow
8+2|z|^2,
$$
where $\lambda_{\max}(\bbH_n)$ denotes the maximum eigenvalue of
$\bbH_n$. It follows that for any $\delta>0$ and sufficiently large $n$
\begin{eqnarray*}
&&\big|\int^\infty_{n^{-4-2\delta}}\log x
\big(v_n(dx,z)-v(dx,z)\big)\big|\\
&=&\big|\int^{8+2|z|^2+\delta}_{n^{-4-2\delta}}\log x\big(v_n(dx,z)-v(dx,z)\big)\big|\\
&\leq& \left(|\log (n^{-4-2\delta})|+\log (8+2|z|^2+\delta)
\right)\|v_n(x,z)-v(x,z)\|\\
 &\stackrel{a.s.}\longrightarrow& 0.
\end{eqnarray*}
Here we do not present the proof of the convergence of $v_n(x,z)$ to
$v(x,z)$ with the desired convergence rate for each $z$. Indeed,
 the rank inequality (see Theorem 11.43 in \cite{bai3})
 can be used to  re-centralize $X_{jk}$
 and then Lemma 10.15 in \cite{bai3} provides
the convergence rate under the assumption $E|X_{11}|^{2+\delta}<\infty$.

On the other hand, by (\ref{a35}) and Borel-Cantelli lemma,
$$
\frac{1}{2n}\log\Big|\det(\bbH_n)\Big|I(s_n(\bbX-z\sqrt{n}\bbI)<n^{-3/2-\delta})\stackrel{a.s.}\longrightarrow0.
$$
Here we take $\varepsilon=n^{-1-\delta},\delta>0$ in (\ref{a35}).
One should observe that $\varepsilon$ in Theorem 5.1 in \cite{rud2}
can be dependent on $n$, so does $\varepsilon$ in Theorem \ref{small}. Moreover,
from Lemma 4.2 in \cite{bai2} one can conclude that
$$
\int^{n^{-4-2\delta}}_0\log x v(dx,z)\rightarrow 0.
$$
Therefore
\begin{equation}
U^{\mu_n(x,y)}\stackrel{a.s.}\longrightarrow
-\frac{1}{2}\int^\infty_0\log x v(dx,z).\label{a30}
\end{equation}
Again by the fourth moment condition
$$
|\lambda_1(\bbX)|\leq \big(\lambda_{\max}(n^{-1}\bbX\bbX^*)\big)^{1/2}
\stackrel{a.s.}\longrightarrow
2.
$$
So for all large $n$, almost surely $\mu_n$ is compactly supported on the disk
$\{z:|z|\leq 2+\delta\}$. Here we have used the fact that all the eigenvalues of
an $n\times n$ matrix are dominated by the largest singular value of the same matrix.
Consequently Theorem \ref{theo3} follows from Lemma \ref{lem4} combined with
Lower Envelop Theorem and Unicity Theorem for logarithmic potential of measures
(see Theorem 6.9, p.73,
and Corollary 2.2, p.98, in \cite{saff}).

\section*{Acknowledgments} The authors would like to thank Prof.
Z. D. Bai for his  helpful discussions when we read Chapter 10 of
Bai and Silverstein's book.

\end{document}